\documentclass[11pt,leqno]{article}
\usepackage{amsmath,amssymb}
\usepackage{epic}
\usepackage{epsfig}

\begin{document}

\newtheorem{defi}{Definition}[section]
\newtheorem{rema}{Remark}[section]
\newtheorem{prop}{Proposition}[section]
\newtheorem{lem}{Lemma}[section]
\newtheorem{theo}{Theorem}[section]
\newtheorem{cor}{Corollary}[section]
\newtheorem{conc}{Conclusion}[section]

\author{P. Germain}

\title{Weak-strong uniqueness for the isentropic compressible Navier-Stokes system}

\maketitle

\maketitle

\begin{abstract}
We prove weak-strong uniqueness results for the isentropic compressible Navier-Stokes system on the torus.
In other words, we give conditions on a strong solution so that it is unique in a class of weak solutions.
Known weak-strong uniqueness results are improved. Classical uniqueness results for this equation follow naturally.
\end{abstract}

\section{Introduction}

\subsection{Presentation of the equation}

We shall study in this paper the Cauchy problem for the isentropic compressible Navier-Stokes system set in the $d$-dimensional torus $\mathbb{T}^d$, with $d\geq 2$. 
The unknowns $\rho$ and $u$ stand for the density and velocity of the fluid; they are respectively $\mathbb{R}^+$, and $\mathbb{R}^d$ valued and they are functions of the space variable $x$ and of the time variable $t$.
The system reads
\begin{equation}
\label{CINS}
\left\{ \begin{array}{l} 
\partial_t \rho + \operatorname{div}( \rho u) = 0 \\
\rho \partial_{t}u + \rho u \cdot \nabla u - L u + \nabla P(\rho) = 0 \\
(\rho,u)_{|t=0} = (\rho_0,u_{0}) \, \, .
\end{array} \right.
\end{equation}
where the notation $u \cdot \nabla$ corresponds to the operator $\sum_{i=1}^d u^i \partial_i$. The pressure law $P$ and the Lam\'e operator $L$ are given by
\begin{equation}
\label{PL}
\begin{split}
& P(\rho) = \rho^\gamma \;\;\;\;\;\mbox{with}\; \gamma > 1 \\
& L = \mu \Delta  + (\lambda + \mu)\nabla \operatorname{div} \;\;\;\;\;\mbox{with}\; \mu>0 \;\;\;\mbox{and}\;\;\; \lambda + 2 \mu > 0 \,\,.
\end{split}
\end{equation}
The pressure depends only on the density due to the isentropy assumption; the Lam\'e operator arises since the fluid is assumed to be Newtonian.

\bigskip

The mass 
$$
\mathcal{M}(\rho) \overset{def}{=} \int \rho \,dx \;\;\;\;\mbox{is constant.}
$$
Formally, energy is also conserved
\begin{equation*}
\mathcal{E}(\rho,u) \overset{def}{=} \frac{1}{2} \|\sqrt{\rho(t)}  u(t)\|_2^2 + \frac{1}{\gamma} \|\rho(t) \|_{\gamma}^{\gamma} + \int_0^t \|\nabla u(s)\|_2^2 \,ds = \frac{1}{2} \|\sqrt{\rho_0} u_0\|_2^2 + \frac{1}{\gamma} \|\rho_0\|_{\gamma}^{\gamma}\,\,.
\end{equation*}
For the weak solutions presented in the next subsection however, it is not known whether the above energy equality holds, or simply the weaker energy inequality, that is the second equality sign above is replaced by $\leq$.

\bigskip

The following transformation leaves the set of solutions invariant:
\begin{equation}
\label{scaling}
\begin{split}
& (u_0,\rho_0) \;\;\mapsto \;\; (\lambda u_0 ( \lambda x), \rho_0 ( \lambda x)) \\
& (u,\rho,P) \;\;\mapsto\;\; (\lambda u (\lambda^2 t , \lambda x), \rho (\lambda^2 t , \lambda x), \lambda^2 P) \,\,.
\end{split}
\end{equation}
This is not a true scaling transformation since the pressure law $P$ has to be scaled too. However if one disregards this problem and only focuses on $u$, one finds the classical scaling of the Navier-Stokes equation. Inspired by this, we will consider that a functional space for $u$ and $\rho$ to be scaling invariant if its norm is left unchanged by the above transformation.

\subsection{Weak and strong solutions}

The first line of research concerning solutions of~(\ref{CINS}) is strong solutions. Though this is no rigorous definition, we mean essentially by this solutions whose uniqueness (in the class where they are built up) can be proved ; they are generally constructed using iteration techniques.

Danchin~\cite{danchin} built up global strong solutions for (small) data at the scaling of the equation. He had to assume that the initial density was bounded by below by a positive number.

Cho, Choe and Kim~\cite{cho} were able to deal with data with vanishing density. However, they required a compatibility condition on their data and could not consider data at the regularity level of the scaling.

\bigskip

Strong solutions have two major shortcomings: their local existence is not known below a certain regularity for the data, or in the presence of vacuum (vanishing density); and their global existence is known only under a smallness assumption.

An explicit example of blow up has even been built up by Weigand~\cite{Weigand}; but one has to apply an unbounded force to the system.

\bigskip

Let us now examine weak solutions. These solutions are in general globally defined, built using a compactness argument, but their uniqueness is not known.

Weak solutions for~(\ref{CINS}) have been proved to exist first under a spherical symmetry assumption by Hoff~\cite{Hoff0}. 

The case of general, finite energy data was treated by Lions~\cite{Lions}, with a technical restriction on $\gamma$, which was reduced by Feireisl~\cite{Feireisl} to $\gamma > \frac{d}{2}$. The weak solution which is obtained has finite energy, but not much more is known of its regularity.

Jiang and Zhang~\cite{Jiang} were able to build up weak solutions for any $\gamma>1$, for axisymmetric data.

Finally, Desjardins~\cite{Desjardins} considered data whose regularity is slightly better than what finite energy would impose, $(\rho_0,u_0) \in L^\infty \times H^1$. He could build up local weak solutions for which this level of regularity is conserved. However, their uniqueness is not known. Notice that $L^\infty$ bound on $\rho_0$ reminds one of the scale-invariant spaces discussed above.

Of course, the basic problem with weak solutions is that the question of their uniqueness is not answered.

\subsection{Weak-strong uniqueness}

The idea of weak-strong uniqueness is the following: assume that a weak solution has some appropriate extra regularity (ie, is also a strong solution in some sense), and prove its uniqueness in the class of finite energy solutions, or some similar weak solution space.
Weak-strong uniqueness gives in particular conditions under which the equation is well behaved, and weak solutions are unique.

The first weak-strong uniqueness results were obtained by Prodi~\cite{Prodi} and Serrin~\cite{Serrin} for the (incompressible) Navier-Stokes equation, for which the same question can be asked. See the articles by the author~\cite{Germain0}~\cite{Germain} for recent results about weak-strong uniqueness for incompressible and inhomogeneous incompressible Navier-Stokes. 

Weak-strong uniqueness criteria for~(\ref{CINS}) were given in~\cite{Desjardins}.

\section{Results obtained}

We present below our two theorems on weak-strong uniqueness for~(\ref{CINS}).

It will be a convenient short-hand to denote
$$
L^p_T L^q \overset{def}{=} L^p([0,T],L^q(\mathbb{T}^d)) \,\,.
$$

\begin{theo}
\label{finitetheo}
Suppose $\gamma > \frac{d}{2}$, and take initial data $(u_0,\rho_0)$ such that
$$
\rho_0 \in L^\gamma \;\;\;\;\mbox{and}\;\;\;\;\sqrt{\rho_0}u_0 \in L^2 \,\,.
$$
A solution $(\bar{u},\bar{\rho})$ is unique on $[0,T]$ in the set of solutions $(\rho,u)$ such that
\begin{equation}
\label{carnaval}
\mathcal{E}(\rho,u) < \infty \;\;\;\mbox{and}\;\;\; \nabla \rho \in L^{2 \gamma}_T L^{\left(\frac{1}{2 \gamma} + \frac{1}{d} \right)^{-1}}
\end{equation}
provided it satisfies, for a constant $c$,  $\bar \rho \geq c > 0$, and
\begin{itemize}
\item if $d = 2$ and $\gamma\geq 2$, $\displaystyle \nabla \bar u \in L^1_T L^\infty \;\;\;\;,\;\;\;\;{L\bar u} \in L^2_T L^p$ with $p>2$.
\item if $d = 2$ and $\gamma< 2$, $\displaystyle \nabla \bar u \in L^1_T L^\infty \;\;\;\;,\;\;\;\;{L\bar u} \in L^2_T L^p \;\;\;\;,\;\;\;\;{L\bar u} \in L^{\frac{2}{3-\gamma}}_T L^q$
with $p>2$, $q>\frac{2}{\gamma-1}$.
\item if $d \geq 3$ and $\gamma\geq 2$, $\displaystyle \nabla \bar u \in L^1_T L^\infty \;\;\;\;,\;\;\;\;{L\bar u} \in L^2_T L^d$.
\item if $d \geq 3$ and $\gamma<2$, $\displaystyle \nabla \bar u \in L^1_T L^\infty \;\;\;\;,\;\;\;\;{L\bar u} \in L^2_T L^d \;\;\;\;,\;\;\;\;{L\bar u} \in L^{\frac{2}{3-\gamma}}_T L^{\frac{d}{\gamma -1}}$.
\end{itemize}
\end{theo}

Remark first that this theorem can be extended to the case where the domain is the whole space $\mathbb{R}^d$ if $d \geq 3$. The proof applies then almost verbatim.

\medskip

Next we would like to mention a remark that was pointed out to us by Nader Masmoudi, and which yields a statement close to the above theorem.

The idea is that the condition $\nabla \rho \in L^{2 \gamma} L^{\left(\frac{1}{2 \gamma} + \frac{1}{d} \right)^{-1}}$ does not actually play any r\^ole in the estimates leading to the theorem, other than justifying formal manipulations.
It therefore should be possible to try and dispense with it. A possibility for doing so is to include the inequality~(\ref{marmotte}), which is the key point in the proof of the theorem, in the definition of a weak solution. To be more specific, we would define a new class of weak solutions, which would not only verify the energy inequality, but also inequality~(\ref{marmotte}) for any $\bar \rho$, $\bar u$ with the appropriate smoothness. 
It would be possible to adapt the argument used to build up global weak solutions in order to include this new condition, since inequality~(\ref{marmotte}) holds for smooth $\rho$, $u$: this is the content of Lemma~\ref{niceform}.

So another version of the above theorem could be obtained, replacing the set given by~(\ref{carnaval}) by the set of solutions $(\rho,u)$ such that
$$
\mathcal{E}(\rho,u) < \infty \;\;\;\mbox{and}\;\;\;(\rho,u) \mbox{ satisfies inequality}~(\ref{marmotte})\,\,. 
$$
A striking fact is that in this new version of Theorem~\ref{finitetheo}, weak strong uniqueness holds without any regularity requirement on the density.

\bigskip

We now come to our second theorem.

\begin{theo}
\label{infinitetheo}
Take initial data $(u_0,\rho_0)$ such that
$$
\rho_0 \in L^\infty \;\;\;\;\mbox{and}\;\;\;\;\sqrt{\rho_0}u_0 \in L^2 \,\,.
$$
A solution $(\bar \rho,\bar u)$ is unique on $[0,T]$ in the set of solutions $(\rho,u)$ such that
$$
\rho \in L^\infty_T L^\infty \;\;\;\; \sqrt{\rho} u \in L^\infty_T L^2 \;\;\;\; \nabla u \in L^2_T L^2 
$$
provided that
\begin{itemize}
\item if $d=2$,
\begin{equation*}
\begin{aligned}
& \nabla \bar \rho \in L^\infty_T L^p \\
& \nabla \bar u \in L^1_T L^\infty \\
& \sqrt{t} \left( \partial_t \bar u + \bar u \cdot \nabla \bar u \right) \in L^2_T L^p \,\,,
\end{aligned}
\end{equation*}
with $p>2$.
\item if $d\geq 3$,
\begin{equation*}
\begin{aligned}
& \nabla \bar \rho \in L^\infty_T L^d \\
& \nabla \bar u \in L^1_T L^\infty \\
& \sqrt{t} \left( \partial_t \bar u + \bar u \cdot \nabla \bar u \right) \in L^2_T L^d \,\,.
\end{aligned}
\end{equation*}
\end{itemize}
\end{theo}

Finally, it is interesting to compare the two theorems above...
\begin{itemize}
\item ... but let us compare them first to existing results. 

Theorem~\ref{finitetheo}, up to the technical condition $\nabla \rho \in L^{2 \gamma}_T L^{\left(\frac{1}{2 \gamma} + \frac{1}{d} \right)^{-1}}$, only requires that the energy of $(u,\rho)$ is finite, which is a very natural condition ; this is the first uniqueness result in such a wide class. Also notice that its proof is quite elementary. There are uniqueness results in classes of smoother functions : Danchin~\cite{danchin} and Hoff~\cite{Hoff}. A nice feature of the result of Hoff is that his criterion includes solutions that exhibit a particular kind of singularity, which is known to exist for $(CINS)$: see Hoff~\cite{Hoff2}.

As for Theorem~\ref{infinitetheo}, it is the first weak-strong uniqueness criterion at the scaling of the equation ; the result of Desjardins~\cite{Desjardins} is close, but misses the scaling.

\item Theorem~\ref{finitetheo} takes advantage of the whole energy (see the proof), whereas Theorem~\ref{infinitetheo} mainly uses the finiteness of the kinetic energy $\|\sqrt{\rho} u \|_2^2$. For this reason, the first theorem above is not restricted to bounded densities. But the exponent $\gamma$ is restricted to $\gamma > \frac{d}{2}$. Notice that this threshhold for $\gamma$ is also the one below which existence of weak solutions is not known.
\item In Theorem~\ref{infinitetheo}, one imposes that $\rho$ is bounded. In some sense, it makes the pseudo-scaling transformation behave like a real scaling transformation (that is, it cancels the trouble with $P$), and endows the equation with a real scaling. It is then no wonder that the conditions provided by Theorem~\ref{infinitetheo} correspond to this scaling (think of the case of Navier-Stokes, where the classical weak-strong uniqueness conditions are at the scaling of the equation).
\item Finally, Theorem~\ref{infinitetheo} includes the possibility of vanishing viscosity $\bar \rho$. The estimates of Theorem~\ref{finitetheo} do not require that $\bar \rho$ stays away from zero, but making the approximation argument work in this case would probably impose heavy additional technical restrictions, so we chose not to consider this case.
\end{itemize}

\section{Proof of Theorem~\ref{finitetheo}}

We consider two solutions of $(CINS)$, $(u,\rho)$ and $(\bar{u},\bar{\rho})$, which share the same initial data, and are both assumed to be of finite energy. In the following, $(u,\rho)$ will play the role of a general finite energy solution, whereas $(\bar{u},\bar{\rho})$ will have additional regularity; our aim will be to prove the uniqueness of $\bar{u}$, that is to prove that $u = \bar{u}$. In order to achieve it, we shall estimate their difference with the help of the following lemma

\begin{lem} 
\label{niceform}
Let $(u,\rho)$ and $(\bar{u},\bar{\rho})$ be two solutions of $(CINS)$ verifying the assumptions of Theorem~(\ref{finitetheo}), and 
$$
U \overset{def}{=} u - \bar{u} \;\;\;\;\;\;\;\;\;\; R \overset{def}{=} \rho - \bar{\rho} \,\,.
$$
Set
$$
F(\bar{\rho},R) = \frac{1}{\gamma} (R + \bar \rho)^\gamma - \bar{\rho}^{\gamma - 1} R - \frac{1}{\gamma} \bar{\rho}^\gamma \,\,.
$$
Then
\begin{equation}
\label{marmotte}
\begin{split}
\partial_t \|\sqrt{\rho} U\|_2^2 + \epsilon \|\nabla U\|_2^2 + \frac{\gamma}{\gamma - 1} & \partial_t \|F(\bar \rho,R)\|_1 \\
& \leq - \int \left( \rho(U \cdot \nabla \bar u) \cdot U + R \frac{L\bar u}{\bar \rho} \cdot U + \gamma \operatorname{div} \bar u F(\bar \rho, R) \right) \,dx \,\,.
\end{split}
\end{equation}
\end{lem}

\begin{rema} 
The quantity $\|\sqrt{\rho} U\|_2^2 + \frac{\gamma}{\gamma - 1}\|F(\bar \rho,R)\|_1$ is sometimes called relative entropy of $\rho$ with respect to $\bar \rho$. It seems to have been first used to prove weak-strong uniqueness results by Dafermos~\cite{Dafermos}, who was studying conservation laws. More recently, Mellet and Vasseur~\cite{MelletVasseur}, by making use of this quantity, obtained weak-strong uniqueness results for the one dimensional isentropic compressible Navier-Stokes equation. This relative entropy is useful in many other contexts : for instance, Berthelin and Vasseur~\cite{BerthelinVasseur} relied on it to prove convergence of kinetic models to the isentropic compressible Euler equation.
\end{rema}

\textsc{Proof of Lemma~\ref{niceform}: }
We will prove the above Lemma without caring about regularity issues when manipulating expressions.
A regularization procedure should be performed in order to prove it with full rigor.
The assumption that $u$ and $\bar u$ satisfy the hypotheses of Theorem~\ref{finitetheo} provides bounds on $u$ and $\bar u$ that would enable one to carry out such a regularization step; we shall however skip this.

First, subtracting the mass conservation equations for $(u,\rho)$ and $(\bar{u},\bar{\rho})$ gives
\begin{equation}
\label{diffmass}
\partial_t R + \operatorname{div}(\rho U) + \operatorname{div}(R \bar u) = 0
\end{equation}
Subtracting the momentum conservation equations gives 
\begin{equation}
\label{butterfly}
(\rho \partial_t + \rho u \cdot \nabla) U - L U + \nabla \rho^\gamma - \nabla \bar{\rho}^\gamma = - R (\partial_t \bar{u} + \bar{u} \cdot \nabla \bar{u}) - \rho U \cdot \nabla \bar{u}\,\,.
\end{equation}
Substituting $\frac{L \bar{u} - \nabla \bar{\rho}^\gamma}{\bar{\rho}}$ for $\partial_t \bar{u} + \bar{u} \cdot \nabla \bar{u}$, this yields
\begin{equation}
\label{coucou1}
(\rho \partial_t + \rho u \cdot \nabla) U - L U + \nabla \rho^\gamma - \frac{\rho}{\bar{\rho}} \nabla \bar{\rho}^\gamma = - \rho U \cdot \nabla \bar{u} - R \frac{L \bar{u}}{\bar{\rho}} \,\,.
\end{equation}
The next step is to take the (space) $L^2$ scalar product of the above with $U$. Using the mass conservation equation, the first term gives
\begin{equation}
\label{coucou2}
\int (\rho \partial_t + \rho u \cdot \nabla) U \cdot U \,dx = \frac{1}{2} \partial_t \| \sqrt{\rho} U \|_2^2 \,\,.
\end{equation}
The second term contributes, due to~(\ref{PL})
\begin{equation}
\label{coucou3}
- \int L U \cdot U\,dx \geq \epsilon \|\nabla U\|_2^2\,\,,
\end{equation}
with $\epsilon >0$. The scalar product of $U$ with the third and fourth summands of~(\ref{coucou1}) is a bit more tricky:
\begin{equation}
\label{coucou4}
\begin{split}
\int & \left( \nabla \rho^\gamma - \frac{\rho}{\bar{\rho}} \nabla \bar{\rho}^\gamma \right) \cdot U \,dx  = \frac{\gamma}{\gamma - 1} \int \rho U \nabla \left( \rho^{\gamma - 1} - \bar{\rho}^{\gamma - 1} \right) \,dx \\
& = - \frac{\gamma}{\gamma - 1} \int \operatorname{div}(\rho U) \left( \rho^{\gamma - 1} - \bar{\rho}^{\gamma - 1} \right) \,dx \;\;\;\;\mbox{integrating by parts}\\
& = \frac{\gamma}{\gamma - 1} \int (\partial_t R + \operatorname{div}(R\bar u))(\rho^{\gamma-1} - \bar{\rho}^{\gamma-1}) \, dx \;\;\;\;\mbox{by~(\ref{diffmass})} \\
& = \frac{\gamma}{\gamma - 1} \left[ \int \partial_t R \frac{\partial}{\partial R} F \,dx + \int \bar u \cdot \nabla R \frac{\partial}{\partial R} F \,dx + \int \operatorname{div} \bar u R \frac{\partial}{\partial R} F \,dx \right] \\
& = \frac{\gamma}{\gamma - 1} \left[ \partial_t \int F \,dx - \int \partial_t \bar \rho \frac{\partial}{\partial \bar \rho} F\,dx + \int \bar u \cdot \nabla F \,dx - \int \bar u \cdot \nabla \bar \rho \frac{\partial}{\partial \bar \rho} F\,dx + \int \operatorname{div} \bar u R \frac{\partial}{\partial R} F \,dx \right] \\
& = \frac{\gamma}{\gamma - 1} \left[ \partial_t \int F \,dx + \int \operatorname{div} \bar u \left( -F + \bar \rho \frac{\partial}{\partial \bar \rho} F + R \frac{\partial}{\partial R} F \right) \,dx \right] \\
& = \frac{\gamma}{\gamma - 1} \partial_t \int F \,dx + \gamma \int \operatorname{div} \bar u F \, dx \,\,.
\end{split}
\end{equation}
Notice that the assumption that $\nabla \rho$ as well as $\nabla \bar\rho$ belong to $L^{2 \gamma}_T  L^{\left(\frac{1}{2 \gamma} + \frac{1}{d} \right)^{-1}}$ implies $\rho$ and $\bar \rho$ belong to $L^{2 \gamma}_T L^{2 \gamma}$. These two facts, combined with $\nabla u \;, \;\nabla \bar u \; \in L^2_T L^2$ ensure that all the expressions written above converge properly.

Now putting together~(\ref{coucou1})-(\ref{coucou2})-(\ref{coucou3})-(\ref{coucou4}), we get the desired estimate. $\blacksquare$

\begin{lem} 
\label{elcomp}
With $F$ defined as in Lemma~\ref{niceform}, for any $\gamma > 1$ there exists a constant $C$ such that
\begin{equation*}
\begin{split}
& R^2 \leq C F(\bar{\rho},R) \bar{\rho}^{2-\gamma} \;\;\;\;\mbox{if} \;R\in[-\bar{\rho},\bar{\rho}] \\
& R^\gamma \leq C F(\bar{\rho},R) \;\;\;\;\mbox{if} \;R \geq \bar{\rho}\,\,.
\end{split}
\end{equation*}
\end{lem}
\textsc{Proof:} 
It only consists of elementary computations, so we skip it. $\blacksquare$

\begin{rema} Since $\rho$ and $\bar \rho$ are non-negative functions, $R = \rho - \bar \rho$ is always larger than $-\bar \rho$.
\end{rema}

As a last ingredient of the proof of Theorem~\ref{finitetheo}, we need the following lemma, almost identical to the one proved in~\cite{Lions1}, page 78:

\begin{lem}
\label{sobolev}
Take $\gamma > \frac{d}{2}$, and suppose that 
$$
\mathcal{M}(\rho) = \int \rho \,dx \geq m \;\;\;\;\mbox{and} \;\;\;\; \int_{\mathbb{R}^d} \rho^\gamma \,dx \leq M
$$
for positive constants $m$ and $M$. Then there exists a constant $C = C(m,M)$ such that
$$
\|v\|_{L^{\frac{2d}{d-2}}} \leq C \left(\|\sqrt{\rho} v\|_2 + \|\nabla v\|_2 \right) \,\,.
$$
\end{lem}
\textsc{Proof:} Let us fix $m$ and $M$ and argue by contradiction. If the lemma is not true, there exists a sequence $(\rho_n,v_n)$ such that
\begin{equation*}
\int \rho^n\,dx \geq m > 0\;\;,\;\; \int \rho_n^\gamma \, dx \leq M \;\;,\;\; \|\sqrt \rho^n v^n \|_2 + \|\nabla v_n\|_2 \rightarrow 0 \;\;\mbox{and}\;\; \|v^n\|_{L^{\frac{2d}{d-2}}} = 1 \,\,.
\end{equation*}
We can extract a subsequence such that, for some $\alpha \in \mathbb{R}$ and $\rho \in L^\gamma$,
\begin{equation*}
\rho^n \rightarrow \rho \;\;\;\;\mbox{weakly in} \; L^\gamma \;\;\;\;\mbox{and}\;\;\;\;\int v_n\,dx \rightarrow \alpha \,\,.
\end{equation*}
This implies
\begin{equation*}
\begin{aligned}
\int \rho \,dx \geq m \;\;\;\;\mbox{and}\;\;\;\;v^n \rightarrow \alpha \;\;\;\;\mbox{strongly in} \; L^{\frac{2d}{d-2}} \,\,.
\end{aligned}
\end{equation*}
So $\alpha$ cannot be zero. Since $\gamma>\frac{d}{2}$ we deduce that
$$
\rho^n |v^n|^2 \rightarrow |\alpha|^2 \rho \;\;\;\;\mbox{weakly in}\;L^1 \,\,,
$$
contradicting the fact that $\|\sqrt \rho^n v^n \|_2 \rightarrow 0$. $\blacksquare$

\bigskip

\textsc{Proof of Theorem~\ref{finitetheo}}: 
We present the proof only in the case where $d \geq 3$. The result is different in dimension $2$ for then the Sobolev embedding $\dot{H}^1 \hookrightarrow L^{\frac{2d}{d-2}}$ does not hold anymore; so one has to modify the proof in an obvious way.

\medskip

We will prove that the right hand side of~(\ref{marmotte}) can be bounded by
$$
\frac{\epsilon}{2} \|\nabla U\|_2^2 + f(t) \left( \|\sqrt{\rho} U\|_2^2 + \|F\|_1 \right) \,\,,
$$
with $f \in L^1([0,T])$; Theorem~(\ref{finitetheo}) follows then by Gronwall's lemma.

\medskip

Bounding the first and third summands of the right-hand side of~(\ref{marmotte}) is fairly easy: 
$$
\int \left( \rho(U \cdot \nabla \bar u) \cdot U + \gamma \operatorname{div} \bar u F(\bar \rho, R) \right) \,dx \leq C \|\nabla \bar{u}\|_{L^\infty} \left( \|\sqrt{\rho} U\|_2^2 + \|F\|_1 \right) \,\,.
$$

\medskip

For the second summand, let us deal first with the region where $R \leq \bar \rho$. One has then, by Lemma~\ref{elcomp}, $R \leq C \sqrt F \bar{\rho}^{1-\gamma/2}$, hence
\begin{equation*}
\begin{split}
\left| \int_{\{R \leq \bar{\rho}\}} R \frac{L\bar u}{\bar \rho} \cdot U \, dx \right| & \leq C \int  \sqrt F \left| \frac{L\bar u}{\bar{\rho}^{\gamma/2}} \right| | U | \, dx \\
& \leq C \|F\|_1^{1/2} \|U\|_{\frac{2d}{d-2}} \left\| \frac{L \bar u}{\bar{\rho}^{\gamma/2}}\right\|_d \;\;\;\;\mbox{by H\"older's inequality} \\
& \leq C \|F\|_1^{1/2} \left( \|\nabla U\|_2 + \|\sqrt \rho U\|_2 \right) \left\| \frac{L \bar u}{\bar{\rho}^{\gamma/2}}\right\|_d \;\;\;\;\mbox{by Lemma~\ref{sobolev}} \\
& \leq C \left\| \frac{L \bar u}{\bar{\rho}^{\gamma/2}}\right\|_d^2 \|F\|_1 + \frac{\epsilon}{4} \left( \|\nabla U\|_2 + \|\sqrt \rho U\|_2 \right)^2 \,\,.
\end{split}
\end{equation*}

\bigskip

Let us now consider the region where $R \geq \bar \rho$. 

If $\gamma \geq 2$, by Lemma~\ref{elcomp}, the inequality $R \leq C \sqrt F \bar{\rho}^{1-\gamma/2}$ still holds, so we find, proceeding as above,
$$
\left| \int_{\{R > \bar{\rho}\}} R \frac{L\bar u}{\bar \rho} \cdot U \, dx \right| \leq C \left\| \frac{L \bar u}{\bar{\rho}^{\gamma/2}}\right\|_d^2 \|F\|_1 +  \frac{\epsilon}{4} \left( \|\nabla U\|_2 + \|\sqrt \rho U\|_2 \right)^2\,\,. 
$$
If $1<\gamma<2$, we notice that $R \leq \rho$ which gives, combined with Lemma~\ref{elcomp}, that $R \leq \sqrt{F} \rho^{1-\gamma/2}$ on $\{R \geq \bar{\rho}\}$. Therefore
\begin{equation*}
\begin{split}
\left| \int_{\{R > \bar{\rho}\}} R \frac{L\bar u}{\bar \rho} \cdot U \, dx \right| & \leq \int \sqrt{F} |U|^{\gamma - 1} \rho^{1-\gamma/2} |U|^{2-\gamma} \left| \frac{L \bar u}{\bar{\rho}} \right| \,dx \\
& \leq C \|F\|^{1/2}_{1} \|U\|_{\frac{2d}{d-2}}^{\gamma-1} \|\sqrt{\rho} U \|_{2}^{2-\gamma} \left\|\frac{L\bar u}{\bar \rho} \right\|_{\frac{d}{\gamma-1}}\;\;\;\;\mbox{by H\"older's inequality} \\
& \leq C \|F\|^{1/2}_{1} \left( \|\nabla U\|_2 + \|\sqrt \rho U\|_2 \right)^{\gamma-1} \|\sqrt{\rho} U \|_{2}^{2-\gamma} \left\|\frac{L\bar u}{\bar \rho} \right\|_{\frac{d}{\gamma-1}}\;\;\;\;\mbox{by Lemma~\ref{sobolev}} \\
& \leq \frac{\epsilon}{4} \left( \|\nabla U\|_2 + \|\sqrt \rho U\|_2 \right)^2 + C \left\| \frac{L\bar u}{\bar \rho} \right\|_{\frac{d}{\gamma-1}}^{\frac{2}{3-\gamma}} \left( \|F\|_1 + \|\sqrt{\rho} U \|_2^2 \right) \,\,.
\end{split}
\end{equation*}

\section{Proof of Theorem~\ref{infinitetheo}}

We shall only consider the case $d \geq 3$; the case $d=2$ can be proved following the same path.

We consider a solution of~(\ref{CINS}), $(\bar{\rho} , \bar u)$, satisfying the assumptions of the Theorem~\ref{infinitetheo}, and would like to prove its uniqueness in the class 
\begin{equation*}
\label{class}
\rho \in L^\infty_T L^\infty \;\;\;\; \sqrt{\rho} u \in L^\infty_T L^2 \;\;\;\; \nabla u \in L^2_T L^2 
\end{equation*}
So let us pick $(\rho,u)$ a solution of~(\ref{CINS}) in the above class, sharing the same data as $(\bar{\rho} , \bar u)$. In order to prove that $(\rho,u) = (\bar{\rho} , \bar u)$, we shall estimate the difference
$$
U \overset{def}{=} u-\bar u \;\;\;\; R \overset{def}{=} \rho - \bar \rho \,\,.
$$

Let us begin by estimating $R$. 

\begin{lem}
\label{controlR}
The following inequality holds
$$
\partial_t \| R \|_2 \leq C \left( ( \|\rho \|_\infty + \| \bar \rho \|_\infty)  \| \nabla U \|_2 + \|\nabla \bar u \|_\infty \|R\|_2 + \|\nabla \bar \rho \|_d \| U \|_{\frac{2d}{d-2}} \right) \,\,.
$$
\end{lem}

\textsc{Proof:} Subtracting the mass conservation equations for $(\rho,u)$ and $(\bar \rho,\bar u)$, we get
$$
\partial_t R + \rho \operatorname{div} U + U \cdot \nabla R + \bar u \cdot \nabla R + R \operatorname{div} \bar u + U \cdot \nabla \bar \rho = 0\,\,.
$$
Taking the (space) scalar product with $R$ and integrating by parts when needed gives the desired result. $\blacksquare$

\bigskip

Let us now estimate $U$.

\begin{lem}
\label{controlU}
If $\rho$ and $\bar \rho$ are bounded by $M$ in $L^\infty$, there exists $C = C(M)$ such that for some $\epsilon >0$
$$
\frac{1}{2}\partial_t \|\sqrt{\rho} U \|_2^2 + \epsilon \|\nabla U\|_2^2 \leq \|R\|_2 \|\partial_t \bar u + \bar u \cdot \nabla \bar u\|_d
\|U\|_{\frac{2d}{d-2}} + C \|\nabla U\|_2 \|R\|_2 + \|\sqrt{\rho} U \|_2^2 \|\nabla \bar u \|_\infty\,\,.
$$
\end{lem}

\textsc{Proof:} This inequality follows in a straightforward fashion after taking the (space) scalar product of~(\ref{butterfly}) with $U$ and noticing that under the assumptions of the lemma,
$$
\left| \langle \nabla \rho^\gamma - \nabla \bar{\rho}^\gamma\,,\,U \rangle \right| \leq C \|R\|_2 \|\nabla U\|_2 \,\,. \,\,\blacksquare
$$

\bigskip

We also need to control the mean value of $U$. This is achieved by the following

\begin{lem}
\label{meanvalue}
The following inequality holds
$$
\left| \int U \,dx \right| \leq \frac{C}{\mathcal{M}(\rho)} \left( \|\rho\|_\infty \|\nabla U\|_2 + \|\nabla \bar u\|_2 \|R\|_2 \right)
$$
\end{lem}

\textsc{Proof:} The lemma follows immediately from the following formula (appearing in Desjardins~\cite{Desjardins})
$$
\int U \,dx = -\frac{1}{\mathcal{M}(\rho)} \int \left( \rho (U - \int U \,dx) + R (\bar u - \int \bar u \,dx) \right) dx \,\,,
$$
and H\"older and Sobolev inequalities. $\blacksquare$

\bigskip

As is always the case for this kind of problem, the proof of the theorem will be concluded by a differential inequality. The classical Gronwall lemma will not suffice here, so we will resort to the more general lemma proved in a previous paper by the author~\cite{Germain}.

\begin{lem}
\label{gronwallplus}
Suppose that the following inequality holds
\begin{equation}
\label{perroquet}
f'+ \left( g'\right)^2 \leq \alpha f + \beta g'g \,\,.
\end{equation}
where $f$, $g'$, $\alpha$, and $\beta$ are positive functions of the real variable such that
$$
f\in L^\infty \;\;\;,\;\;\;g(0)=0\;\;\;,\;\;\; g' \in L^2\;\;\;,\;\;\;\alpha \in L^1 \;\;\;\mbox{and}\;\;\;\sqrt{t}\beta(t)\in L^2\,\,.
$$
Then there holds
\begin{equation*}
\label{cacatoes}
\left( e^{-\int_0^t \alpha} - \frac{1}{2}-\frac{1}{2}\int_0^t s \beta(s)^2 \,ds \right)
\int_0^t \left( g'\right)^2 + e^{-\int_0^t \alpha} f(t) \leq f(0) \,\,.
\end{equation*}
\end{lem}

\bigskip

We can now proceed with the proof of the theorem.

\bigskip

\textsc{Proof of Theorem~\ref{infinitetheo}:} We shall from now on assume that $(\rho,u)$ and $(\bar \rho , \bar u)$ satisfy the assumptions of the theorem.

First, Lemma~\ref{meanvalue} and Sobolev inequality yield the control of $\|U\|_{\frac{2d}{d-2}}$:
\begin{equation}
\label{beetle}
\|U\|_{\frac{2d}{d-2}} \leq C \left( \|\nabla U\|_2 + \|\nabla \bar u\|_2 \|R\|_2 \right) \,\,.
\end{equation}
Inserting this in Lemma~\ref{controlR} gives
$$
\partial_t \|R\|_2 \leq C \left( \| \nabla U \|_2 + \left( \|\nabla \bar u\|_\infty + \|\nabla \bar u\|_2 \right) \|R\|_2 \right) \,\,.
$$
Integrating this last inequality, we obtain
$$
\|R\|_2(t) \leq C \int_0^t \| \nabla U \|_2 \,ds \,\,.
$$
If we use this last inequality to estimate $\|R\|_2$ and~(\ref{beetle}) to control $\|U\|_{\frac{2d}{d-2}}$, Lemma~\ref{controlU} becomes
\begin{equation*}
\begin{split}
\frac{1}{2} \partial_t & \|\sqrt \rho U \|_2^2  + \epsilon \|\nabla U \|_2^2 \\
& \leq C \left( \|\partial_t \bar u + \bar u \cdot \nabla \bar u\|_d \| \nabla U\|_2 \int_0^t \| \nabla U \|_2 \,ds + \|\partial_t \bar u + \bar u \cdot \nabla \bar u\|_d \|\nabla \bar u\|_2 \left( \int_0^t \| \nabla U \|_2 \,ds \right)^2 \right. \\
&\;\;\;\;\;\;\;\; + \left. \| \nabla U\|_2 \int_0^t \| \nabla U \|_2 \,ds + \| \nabla \bar u \|_{\infty} \| \sqrt \rho U \|_2^2 \right) \\
& \leq C \left( (\|\partial_t \bar u + \bar u \cdot \nabla \bar u\|_d + 1) \| \nabla U\|_2 \int_0^t \| \nabla U \|_2 \,ds + t \|\partial_t \bar u + \bar u \cdot \nabla \bar u\|_d \|\nabla \bar u\|_2 \left( \int_0^t \|\nabla U\|_2^2 \,ds \right) \right.\\
&\;\;\;\;\;\;\;\; \left. + \| \nabla \bar u \|_{\infty} \| \sqrt \rho U \|_2^2 \right) \,\,.
\end{split}
\end{equation*}
In order to conclude, it suffices to apply Lemma~\ref{gronwallplus} with
\begin{equation*}
\begin{aligned}
& f(t) = \frac{1}{2}\|\sqrt \rho U (t) \|_2^2 + \frac{\epsilon}{2} \int_0^t \|\nabla U\|_2^2 \,ds \\
& g'(t) = \frac{\sqrt{\epsilon}}{\sqrt{2}}\|\nabla U\|_2 \\
& \alpha = 2 \| \nabla \bar u\|_\infty + \frac{2}{\epsilon} t \|\nabla \bar u\|_2 \|\partial_t \bar u + \bar u \cdot \nabla \bar u\|_d \\
& \beta = \frac{2}{\epsilon} \left( \| \partial_t \bar u + \bar u \cdot \nabla \bar u \|_d + 1 \right) \,\,.\,\,\blacksquare
\end{aligned}
\end{equation*}

\bigskip

\bigskip

{\bf Acknowledgement:} The author wishes to thank Nader Masmoudi for very helpful comments during the writing of this paper.

\bigskip

\bigskip

\bigskip

Pierre GERMAIN

{\sc Courant Institute of Mathematical Sciences

New York University 

New York, NY 10012-1185

USA}

\bigskip
  
{\tt pgermain@math.nyu.edu}

\end{document}